# A COGNITIVE ANALYSIS OF CAUCHY'S CONCEPTIONS OF FUNCTION, CONTINUITY, LIMIT, AND INFINITESIMAL, WITH IMPLICATIONS FOR TEACHING THE CALCULUS

David Tall

Mathematics Education Research Centre
University of Warwick
CV4 7AL, United Kingdom
<david.tall@warwick.ac.uk>

Mikhail Katz

Department of Mathematics,
Bar Ilan University,
Ramat Gan 52900 Israel
<katzmik@macs.biu.ac.il>

*In this paper we use theoretical frameworks from mathematics education and cognitive psychology to analyse Cauchy's ideas of function, continuity, limit and infinitesimal expressed in his Cours D'Analyse. Our analysis focuses on the development of mathematical thinking from human perception and action into more sophisticated forms of reasoning and proof, offering different insights from those afforded by historical or mathematical analyses. It highlights the conceptual power of Cauchy's vision and the fundamental change involved in passing from the dynamic variability of the calculus to the modern set-theoretic formulation of mathematical analysis. This offers a re-evaluation of the relationship between the natural geometry and algebra of elementary calculus that continues to be used in applied mathematics, and the formal set theory of mathematical analysis that develops in pure mathematics and evolves into the logical development of non-standard analysis using infinitesimal concepts. It suggests that educational theories developed to evaluate student learning are themselves based on the conceptions of the experts who formulate them. It encourages us to reflect on the principles that we use to analyse the developing mathematical thinking of students, and to make an effort to understand the rationale of differing theoretical viewpoints.*

## 1. INTRODUCTION

The calculus today is viewed in two fundamentally different ways: the theoretical calculus used in applications, often based on the symbolism of Leibniz, enhanced by the work of Cauchy using infinitesimal techniques, and the formal mathematical analysis of Weierstrass based on quantified set theory. In this paper we begin in section 2 by reviewing current scholarship in the history of mathematics. In section 3 we complement the historical approach with an analysis of the evolution of mathematical ideas based on human cognitive development. This begins with the theory of Merlin Donald (2001) who distinguishes three levels of thinking, from immediate perception, to dynamic change over short periods of time and on to extended awareness appropriate for the building of more sophisticated theories. It continues with a review of other theories, including Tall's (2013) framework of three distinct ways of thinking mathematically through embodiment, symbolism and formalism that is consistent with Donald's theory and together they provide a fundamental framework to consider the evolution of ideas from ancient Greece, through the development of the calculus, on to the formal development of mathematical analysis and beyond.

Section 4 offers a detailed analysis of Cauchy's conceptions, arising from the work of his predecessors to build his own view of quantities that could be constant or variable and marked on a number line, allowing him to develop his own conceptions of functions, continuity and limits, using the idea of an infinitesimal generated by a variable sequence whose terms tended to zero.

Section 5 follows the transition to the mathematical analysis of Weierstrass, the formalist approach of Hilbert and the rejection of the infinitesimal approach by Cantor as new approaches evolve in the twentieth century. However, we show that the formalist approach of Hilbert can be used to prove a simple structure theorem that any ordered field containing the real numbers as a subfield must contain infinitesimals and that any element that lies between two real numbers must be of the form 'a real number plus an infinitesimal'. This simple theorem should not be confused with the more complex issues of the logical approach to non-standard analysis.

In section 6, we note that the idea of arbitrarily small quantities that evolved in history continue to arise in our students and to be used in applications of mathematics.

This leads to the final section in which we consider the implications of these results to mathematics education, which suggest that different communities of practice, in school, in various applications, and in university pure mathematics can each hold a particular view of calculus that may be appropriate in their own community but may not be appropriate in others.

## 2. SCHOLARSHIP IN THE HISTORY OF MATHEMATICS

An examination of contemporary Cauchy scholarship reveals a surprising lack of consensus among scholars concerning his contribution to the development of modern views of calculus and analysis. The traditional view elaborated by Grabiner (1983, p. 185), is that

> Delta-epsilon proofs are first found in the works of Augustin-Louis Cauchy.

Schubring (2005, p. 480), on the other hand, based on a study of Cauchy's published works and letters, concludes that

> he was actually confined by his adherence to the Leibnizian-Newtonian tradition of the eighteenth century.

A recent study by Ehrlich (2006) documents an uninterrupted tradition of work on infinitesimal-enriched systems, from the end of the 19$^{th}$ century onward. Meanwhile, traditional accounts tend to credit Cantor, Dedekind, and Weierstrass with the 'elimination' of infinitesimals from mathematics. Cantor, on the one hand, introduced infinite cardinals that did not have multiplicative inverses, and on the other, completed the number line using Cauchy sequences of rationals (Cantor, 1872) to postulate a one-to-one correspondence between real numbers and the linear continuum of geometry. Cantor is believed by many professional mathematicians to have shown that there is only one consistent conception of the geometric continuum, with irrationals 'filling' the gaps between the rationals to 'complete' the real number line and leave no room for ephemeral entities that are 'infinitesimal'. The real numbers of Cantor and Dedekind were accepted as the foundation of modern analysis developed by Weierstrass.

However, Felix Klein reflected on this development and observed that

> The scientific mathematics of today is built upon the series of developments which we have been outlining. But an essentially different conception of infinitesimal calculus has been running parallel with this through the centuries. (Klein, 1908, p. 214.)

He noted that this parallel conception of calculus 'harks back to old metaphysical speculations concerning the structure of the continuum according to which this was made up of ultimate indivisible infinitely small parts.' (ibid. p. 214.)

When Robinson (1966) provided a logical approach to the calculus using infinitesimals, the role of Cauchy's ideas was again brought to the fore. Was Cauchy a forerunner of the epsilon-delta approach of Weierstrass, or of the infinitesimal approach of Robinson? Or was he something else?

In considering various modern views by historians and mathematicians, Grattan-Guinness expanded on an observation by Freudenthal (1971, p. 377), to the effect that succeeding generations interpret earlier mathematics in different ways, in the following terms:

> It is mere feedback-style ahistory to read Cauchy (and contemporaries such as Bernard Bolzano) as if they had read Weierstrass already. On the contrary, their own pre-Weierstrassian muddles need historical reconstruction. (Grattan-Guinness, 2004, p. 176.)

Bair et al. (2013) analyse the assumptions underlying post-Weierstrassian approaches to the history of mathematics and argue for a more sympathetic account of the history of infinitesimal mathematics.

Modern views of Cauchy's work are inevitably seen through the eyes of various commentators each of whom has a cultural background that adopts a particular view of mathematics. For instance, in looking at a geometrical line, do we see a continuum with length and no breadth as originally conceived in Euclid Book I—where points lie *on* lines, either as endpoints or where two lines meet—or do we see the line as a set-theoretic aggregate of the points that make it up? According to Cantor (1872), the real line is precisely the complete ordered field that is the basis of modern epsilon-delta analysis and does not contain any infinitesimals. Yet, according to the non-standard analysis of Robinson (1966), the real line $\mathbb{R}$ is part of an enhanced continuum $^*\mathbb{R}$ which includes not only infinitesimals, but also infinite elements that are their inverses.

How then do we interpret the past taking explicit account of the colouring of our interpretations with the encrustations of later generations? This proves to be a complex and subtle enterprise where contributions may be made by historians, philosophers, mathematicians, cognitive scientists and others, each offering subtly different views with a variety of possibly conflicting interpretations.

## 3. THEORIES OF COGNITIVE DEVELOPMENT

Our approach here is to consider the underlying nature of human perception and how our unique human facility for language and symbolism leads to more sophisticated forms of mathematical reasoning, both in the individual and in shared communities of practice over the generations.

For example, although we may speak about a point in space having a position but no size, none of us can actually *see* such a point because it is too small and we represent it as a physical mark with a finite size. Not only can we not perceive arbitrarily small quantities in space, we cannot perceive arbitrarily small intervals in time because the brain takes around a fortieth of a second to build up a selective binding of neural structures to interpret our perceptions (the precise time depending on the extent of the neural connections involved). This leads to a conflict between our perception, in which a halving process soon produces a perceptual quantity that is too small to see, and a symbolic process of halving that involves a symbolic quantity that never precisely equals zero. We will consider the underlying mechanisms that shape human thought, building from the foundational level of human perception through successive levels of consciousness identified by the psychologist Merlin Donald.

### 3.1 Donald's three levels of consciousness

In his book *A Mind So Rare*, Donald (2001) suggests that human consciousness works on three levels. The first is the immediate consciousness of our perception, which takes about a fortieth of a second to combine neural activity into a thinkable concept. Indeed, we have known since the development of the moving pictures in the early twentieth century, that showing a sequence of still images separated by a short interval of darkness offers a sense of relatively smooth motion at speeds faster than around fourteen frames per second (the minimum speed in early digital cameras offering moderately smooth video recording).

Our perception of motion occurs by coordinating our changing perception through a second level of consciousness that Donald calls 'short-term awareness', lasting continuously over a period of two or three seconds. It is this level of consciousness that allows us to perceive a shape moving in space, so as to recognise that it is the same entity seen from different viewpoints.

His third level of consciousness, which he terms 'extended awareness', occurs over longer periods as we reflect on our previous experiences, bringing disparate ideas together, possibly through recording them using words, symbols, or pictures that we can consider simultaneously to build more subtle relationships.

The perceptual notion of continuity is a second level operation. Perceptual continuity occurs both as a phenomenon in time, as we dynamically draw a curve with a stroke of a pencil, without taking the pencil off the paper, and also in space as we observe a continuous curve without any gaps, or a continuous surface or solid.

In a similar vein, Klein notes that

> It is precisely in the discovery and in the development of the infinitesimal calculus that this inductive process, built up without compelling logical steps, played such a great role; and the effective heuristic aid was very often sense perception. And I mean here immediate sense perception, with all its inexactness, for which a curve is a stroke of definite width, rather than an abstract perception which postulates a complete passage to the limit, yielding a one dimensional line. (Klein, 1908, p. 208.)

As mathematicians, we build on the perceptual continuity experienced through our short-term awareness, to shift to a third level of extended awareness where formal continuity is defined verbally and symbolically. Such a shift changes the focus on continuity, from a global phenomenon drawing a graph dynamically with a stroke of a pencil, to a formal definition in terms of a Weierstrassian challenge 'tell me how close you want $f(x)$ to be to $f(a)$, and I will tell you how small you need to make the difference between $x$ and $a$ so that $f(x)$ and $f(a)$ are as close as desired.' Such a computable, yet cumbersome, description of continuity has the effect of shifting the mathematics from a perceptual idea related to our senses, to a computational process that is potentially achievable. In the late 19th century, this approach became the foundation of modern mathematical analysis.

### 3.2 Theories of the cognitive development of mathematical thinking

There are various theoretical perspectives that are relevant in considering the development from human perception of change to its symbolisation in calculus and formalisation in mathematical analysis.

Recent theories of the development of mathematical thinking (Dubinsky, 1991; Sfard, 1991; Gray & Tall 1994) focus on the way that humans carry out mathematical operations, such as addition, sharing, calculating the limit of a sequence, differentiation and integration. At each stage, we perform a process that occurs in time, to produce an output that may also be conceived as a mental entity, independent of time. Counting gives rise to the concept of number, the process of addition 3+2 gives rise to a mental concept, the sum, which is also written using the same symbol 3+2; sharing gives the concept of fraction; calculating a trigonometric ratio such as $\sin A$ = opposite/hypotenuse gives the concept of sine; the process of differentiation gives the derivative; the process of integration gives the integral. In every case, a symbolic notation, such as 3+2, ¾, $a(b+c)$, $dy/dx$, $\int f(x)dx$, $\sum_{n=1}^{\infty} 1/n^2$, represents both a desired process and the resulting concept. Gray & Tall (1994) refer to this conception of a symbol, that dually represents process or concept, as a *procept*, with the additional flexiblity that different symbols with the same output, such as $a(b+c)$ and $ab+ac$, represent the same procept.

Analyzing a sequence $(a_n)$ given by a specific formula involves computing a succession of terms and observing how the process tends to a specific value—the limit. Cornu (1983, 1991) described how students think of the process of getting small as producing an object that is itself arbitrarily small, but not zero. Tall (1986) defined this to be a *generic limit*. It occurs when the focus of attention shifts from a sequence of distinct real constants $a_1, a_2, a_3, \ldots, a_n, \ldots$ to considering the term $a_n$ *as a variable quantity* that varies with *n*, and, as *n* increases, the variable quantity becomes smaller but may never become zero. This produces a mental compression of thought in which distinct terms are conceived and spoken of as a single varying entity.

Such a compression of thought is a widespread phenomenon (Cornu, 1991). For instance, the infinite decimal 0.999… is widely considered to be a quantity 'just less than 1', or even 'as close as it is possible to reach 1 without actually equalling it.' This phenomenon has been analysed at a number of levels over the years from the conflict caused by a potentially infinite process that goes on forever to more recent analyses involving the meaning of infinitesimal concepts (eg. Katz & Katz, 2010a, 2010b and Ely, 2010).

Tall (1986) linked the phenomenon to the historic 'principle of continuity' (or 'law of continuity') formulated by Leibniz in an unpublished text referred to as *Cum Prodiisset*, circ. 1701 in the following terms:

> In any supposed [continuous][1] transition, ending in any terminus, it is permissible to institute a general reasoning, in which the final terminus may also be included. (Child, 1920, p. 147.)

A related formulation is found in a 1702 letter to Varignon, where Leibniz wrote:

> The rules of the finite succeed in the infinite […] and vice versa, the rules of the infinite succeed in the finite.

(see Katz and Sherry (2012 and 2013) and Sherry and Katz (2013) for a fuller discussion). Based on this general heuristic, Leibniz was able to assume that infinitesimals enjoy the same properties as ordinary numbers, and to operate on them accordingly.

Lakoff and Núñez reformulated a related principle in their 'basic metaphor of infinity' in which

> We hypothesize that all cases of infinity—infinite sets, points at infinity, limits of infinite series, infinite intersections, least upper bounds—are special cases of a single conceptual metaphor in which processes that go on indefinitely are conceptualized as having an end and an ultimate result. (Lakoff and Núñez, 2000, p. 258.)

The underlying brain activity is more fundamental. Useful links between neurons are strengthened and provide new and more immediate paths of thought in which an on-going process may be imagined as a thinkable concept. Throughout the history of the calculus, mathematicians ranging from Leibniz to Euler to Cauchy sought to motivate infinitesimals in terms of sequences tending to zero. Thus, a potentially infinite temporal sequence of distinct terms ($a_n$) tending to zero is verbalised as a single term $a_n$ that *varies* as *n* increases, producing a mental entity that is arbitrarily small yet not zero.

---

[1] Child omitted the word 'continuous', which is in the original Latin.

This is a conception that has caused controversy for millennia. If one continually halves a quantity, does the process continue forever becoming smaller and smaller without ever reaching zero, or does one reach an indivisible element that can no longer be further subdivided? The perceptual brain, with its finite capacity to perceive and sense, cannot resolve this problem but the brain's extended level of awareness (Donald's level 3) is able to formulate the problem conceptually and linguistically to express opinions and beliefs about a possible resolution.

### 3.3 Tall's three worlds of mathematics

The link from the perceptual world of our human experience to the computational world of arithmetic and algebra was initiated by Vieta, Descartes, and others. This led to the calculus of Newton and Leibniz that computed naturally perceived phenomena such as the measurement of length, area, volume, time, distance, velocity, acceleration and their rate of change and growth using the computational and manipulable symbolism of the calculus. Cantor, Dedekind and Weierstrass took matters a step further by interpreting the number line in terms of the symbolism of number, arithmetic, order and completeness.

However, it was the introduction of formal axiomatic mathematics by Hilbert (1900) that radically changed the way in which we are able to think of mathematical ideas. He switched attention from the natural phenomena that we perceive physically and conceive mentally, to the *properties* of the phenomena. A mathematical structure is specified by axioms, and deductions are made by mathematical proof. This releases mathematical thinking from the limitations of human perception to the possibilities of formally-defined systems and their consequent properties.

Tall (2004, 2008) formulated three essentially different ways in which mathematical thinking develops, which relate both to the historical development of ideas and also to the cognitive development of the individual from child to mathematician:

(1) *Conceptual embodiment* builds on human perceptions and actions, developing mental images that are verbalized in increasingly sophisticated ways and become perfect mental entities in our imagination.

(2) *Proceptual[2] symbolism* grows out of physical actions into mathematical procedures that are symbolized and conceived dually as operations to perform and symbols that can themselves be operated on by calculation and manipulation.

---

[2] The name of the mental world of *proceptual* symbolism was later modified to *operational* symbolism (e.g. in Tall, 2010, 2013) to include the wide range of flexible proceptual and routine procedural operations used by different learners. Here it is entirely proper to focus on the more subtle flexible relationship between process and concept.

(3) *Axiomatic formalism* builds formal knowledge in axiomatic systems in a suitable foundational framework (such as formal set theory or formal logic) whose properties are deduced by mathematical proof.

In what follows we will use the shortened terms, 'embodied', 'symbolic' and 'formal' on the understanding that these words reflect the meanings given here rather than the various other meanings found in the literature. In broad terms, these ways of mathematical thinking develop one from another in the cognitive growth of the child into the adult.

Mathematical ideas begin in the perception and action of the embodied world. Perceptions can develop through description, definition and deduction in Euclidean geometry (as noted in successive levels by van Hiele, 1986), within an increasingly sophisticated world of conceptual embodiment, refined through subtle thought experiment and the use of language to develop platonic conceptions of abstract thought.

Actions can lead to mathematical operations that are symbolised and give new forms of calculation and manipulation at increasingly sophisticated levels. Conceptual embodiment and proceptual symbolism develop in parallel, blending together various aspects to give even more sophisticated forms of mathematical thinking.

Distinct types of mathematical reasoning arise in embodiment and symbolism. Euclidean proof in geometry is expressed verbally in terms of definitions and theorems and yet is fundamentally based on human embodiment such as placing one triangle on top of another to see if they fit, as in the principle of congruence. Proof in arithmetic begins with the observation of regularities such as the properties of multiplication of whole numbers leading to the concepts of prime and composite numbers and the proof that any whole number has a unique expression as a product of primes. Algebraic proof builds on the observed regularities of arithmetical operations to base proofs on the 'rules of arithmetic'.

At a later stage, linguistic and logical sophistication lead to a form of mathematics that is presented in terms of set-theoretic axiomatic definitions and formal proof.

In this framework, the perceptual conception of continuity is now seen in terms of embodiment. It relates to the dynamic movement of a pencil as one draws the graph to get a curve that is perceptually continuous without any gaps. Such an interpretation is often seen as an 'intuitive' form of continuity that lacks formal precision. However, such perceptual beginnings are an essential starting point for the later development of more sophisticated forms of mathematical thinking that develop both in the individual and also in succeeding generations over the centuries.

Cauchy flourished at a time in history when the calculus had developed into a remarkable facility in computing the symbolic solutions to many problems involving change and growth, yet was still subject to controversy over various aspects of the theory, in particular, the use of infinitesimal quantities that could be imagined in the mind yet not perceived by the physical human eye.

In the introduction to his *Cours D'Analyse* of 1821, he declared that he could not speak about the continuity of functions without using the properties of infinitely small quantities:

> En parlant de la continuité des fonctions, je n'ai pu me dispenser de faire connaître les propriétés principales des quantités infiniment petites, propriétés qui servent de base au calcul infinitésimal. (Cauchy, 1821, Introduction p. ii.)[3]

The question is: what *are* Cauchy's 'quantités infiniment petites'? Where do they fit in the development from perceptual embodiment to symbolic process and concept, and on to appropriate forms of mathematical definition and proof? As we reflect on this, we must attempt to consider Cauchy's ideas in his own terms and not simply analyse them in a manner that reflects our current views of mathematics that were not available in his time.

## 4. CAUCHY

Cauchy's ideas arose at a particular time in history as French mathematics benefited from the developments of previous generations, in particular, the insight of Descartes, who united the geometry of the Greeks with the vision of points on the line and in the plane corresponding to numbers. His development should be seen in terms of his context at the time.

### 4.1 The legacy of Cauchy's predecessors

Geometric ideas relate back to the theory of Euclid using definitions and explicit assumptions (called 'postulates' and 'common notions') to build a proven system of propositions deduced from first principles. Implicitly, Euclidean proof depends on being able to interpret the meanings of the definitions. For instance, the first four definitions of Euclid Book I specify a point as 'that which has no part', a line is 'breadthless length', 'the ends of a line are points', and 'a straight line is a line which lies evenly with the points on itself.'

Despite the apparent dependence on the implicit meanings of the terms concerned, these definitions carry a number of implications. For instance, a 'line'—which may be curved—has a (finite) length with a point at each end. A point can lie on a line, (as an endpoint, or a point where lines cross), but a line is an entity in itself and does not 'consist of points' in a modern set-theoretic sense. From a modern viewpoint where a line is seen as consisting of a set of points, this simple idea leads to controversy. For instance, Ferraro (2004, p. 37) argued that a Euclidean geometric line *has* endpoints and so one cannot think of it separately as an open interval without endpoints. On the other hand, Cauchy scholars (such as Grabiner, 1983, Laugwitz, 1987, p. 261) see Cauchy variables as ranging through values represented as variable points on a line without including the endpoints, and argue that both open and closed intervals are found in Cauchy. A level 2 perception of drawing a finite line segment in Euclidean geometry sees a line as an entity in itself drawn from one point to another, and points lie on it, either at the ends or where one line crosses another.

---

[3] The page numbers in *Cours D'Analyse* refer to the collected edition of Cauchy's works that can be downloaded for personal use from http://gallica.bnf.fr/ark:/12148/bpt6k90195m/f12

For the Greeks, geometry and arithmetic were considered as being fundamentally distinct. Arithmetic operations could be applied to whole numbers and fractions, but the product of two lengths is an area, so the Greeks did not build a full arithmetic in geometry. Instead, they compared magnitudes of the same kind (length with length or area with area) as in the following translation presented by Joyce (1997):

> Triangles and parallelograms which have the same height are to one another as their bases. (Euclid Book VI, proposition 1.)

The Greek theory of proportions focused on comparing magnitudes by various methods, for example, using parallel lines in theorems such as:

> If a straight line is drawn parallel to one of the sides of a triangle, then it cuts the sides of the triangle proportionally. (Euclid Book VI, Proposition 2.)

In figure 1 the line DE is constructed to be parallel to the base BC, and the equality of proportions is written as AB:AD :: AC:AE and spoken as 'AB is to AD as AC is to AE'.

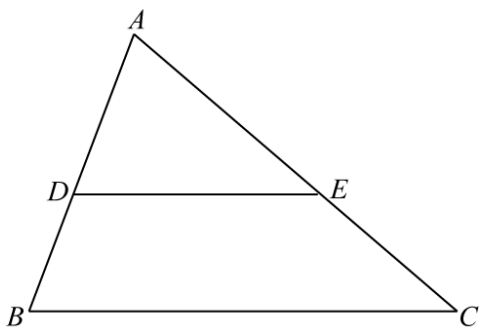

Figure 1: parallel lines and proportion

In the seventeenth century, Descartes (1635) set himself free of the limitations imposed by this view by the simple expedient of choosing a length to be the unit in a particular calculation. He was then able to multiply and divide lengths.

By rewriting AB:AD :: AC:AE as

$$\frac{AB}{AD} = \frac{AC}{AE}$$

he could calculate $AB = AC \times AD$ by choosing $AE = 1$ and performing a suitable geometric construction. Not only did this give a full arithmetic of numbers as lengths, it allowed him to link the position of a point on a line to a signed number using a chosen origin and unit length. This extended to the identification of points in the plane as ordered pairs of numbers relative to a given pair of axes. Curves in the plane could now be described by algebraic relationships, allowing a complementary interchange of ideas between geometry and algebra.

Leibniz (1684) used the ideas of Descartes to compute the slopes of curves and areas that had curves as boundaries, by imagining a curve as a polygon with an infinite number of infinitesimally small sides, so that he could compute the slope of the curve as the quotient *dy/dx* for infinitesimal increments *dx* and *dy* in the independent variable *x* and the dependent variable *y*.

Subsequently Euler (1748) continued the development of calculus in the eighteenth century, focussing on algebra rather than geometry, declaring the principle of 'the generality of algebra' that calculations with complex numbers and infinite series are subject to the same principles as the natural properties of operations with ordinary numbers. Such a principle led him to operate with infinitesimals in a purely symbolic manner, where infinitesimal quantities were conceived to be zero, in the sense that they did not change the value of a finite number when added, and yet operated in such a way that their ratios *dy/dx* could be calculated. His concept of the generality of algebra took him much further, declaring that

> Calculus concerns variable quantities, that is, quantities considered in general. If it were not generally true that d(log *x*)/*dx* = 1/*x*, whatever value we give to *x*, either positive, negative, or even imaginary, we would never be able to make use of this rule, the truth of the differential calculus being founded on the generality of the rules it contains. (Euler, 1749.)

However, such operations with infinitesimals did not meet with universal approval. For instance, d'Alembert observed that, contrary to generally held opinion,

> In the differential calculus one is not at all concerned with infinitely small quantities but only with the limits of definite quantities. The words infinite and small are used only to abbreviate expressions. (Breuss, 1968.)

At the beginning of the nineteenth century, Cauchy adopted a view rooted in variable quantities (already found in l'Hôpital, 1696), by defining an infinitesimal in terms of a sequence of familiar numbers whose values diminished to zero.

**4.2 Cauchy's distinction between number and quantity**

Cauchy lived in a culture that used the Cartesian correspondence between points on a line and (signed) numbers, which extended to provide a correspondence between points in the plane and pairs of (signed) numbers. In this correspondence, points lay geometrically *on* the line or in the plane, the line and the plane were not solely an aggregate of the set of points. Cauchy used these correspondences in calculus and also in the more subtle development of complex differentiation and integration. His concept of magnitude related to fundamental Euclidean meanings where angles, lengths, areas and volumes are magnitudes that are numbers without a sign, but signed numbers could be used as quantities in calculations and to mark the positions of points on a line or in the plane.

In the opening pages of his preliminary chapter of *Cours D'Analyse* (1821), written for his students at the École Polytechnique, he steadily unfolds his distinction between *numbers* which are unsigned and used to count and measure, and *quantities*, which occur in operations with numbers. 'Quantities' may be positive or negative, written as numbers preceded by a sign (+ or –) and marked on the number line—with positive numbers to the right of the origin and negative to the left—operating with prescribed rules, such as 'the product of two negative quantities is positive.'

He then goes on to speak of *constant quantities* and *variable quantities*, where constants have fixed values and are represented by letters at the beginning of the alphabet, and variable quantities can take on various values and are represented by letters at the end of the alphabet. These are long-established conventions that go back to Euler. However, what is distinctive in his work is the way that he uses 'infinitesimal quantities' as *computational* variables.

Cauchy noted that certain variable quantities take on successive values that approach a fixed value (lorsque les valeurs successivement attribuées à une même variable s'approchent indéfiniment d'une valeur fixe), so that they differ from that value by as little as one wishes (de manière à finir par en différer aussi peu que l'on voudra). The fixed value is then called the *limit* of all the others (cette dernière est appelée la *limite* des toute les autres).

Later, in the opening of chapter II, he confirms that his notion of 'approaching indefinitely' (approchent indéfiniment) does not necessarily imply a monotonic sequence where terms are successively smaller, but could include a succession of numbers such as

$$\tfrac{1}{4}, \tfrac{1}{3}, \tfrac{1}{6}, \tfrac{1}{5}, \tfrac{1}{8}, \tfrac{1}{7}, \ldots$$

which eventually become smaller than any given value.

This concept of limit reads very much like our modern notion, which is supported by his example that an irrational number is a limit of a sequence of fractions. But then he gives a further example:

> En Géométrie, la surface du cercle est la limite vers laquelle convergent les surfaces des polygones inscrits, tandis que le nombre de leurs côtés croît de plus en plus.

'A circle is the limit of inscribed polygons as the number of sides increases more and more.' Thus Cauchy's notion of limit already has aspects that extend the idea of a limit of a sequence of numbers to more general mathematical situations. In particular, in this case the polygonal figure grows outwards to get as close as desired to the limiting circle. This particular example has the limiting circle as a barrier that 'limits' the growth of the polygon. In general, Cauchy allows a limit to be approached both from above and below (as in case of a series with terms having alternate signs) but he also distinguishes between the limit of a function, such as $1/x$ at the origin, which he sees as being positive infinite from the left and negative infinite from the right, an example to which we will return later.

At the beginning of his *Résumé des Leçons…*, he introduces the concept of infinitesimal quantity (infiniment petit) in terms of a variable with a sequence of values whose absolute value (valeur numérique) tends to zero:

> Lorsque les valeurs numériques successives d'une même variable décroissent indéfiniment, de manière à s'abaisser au-dessous de tout nombre donné, cette variable

> devient ce qu'on nomme un infiniment petit ou une quantité infiniment petite. Une variablé de cette espèce a zéro pour limite. (Cauchy, 1823, p. 4.)

In this definition, he uses both the term 'infiniment petit' as a *noun*, and 'une quantité infiniment petite' where the phrase 'infiniment petite' is an *adjective* applying to the feminine noun 'quantité'. Conveniently, the gender given by the ending of 'petit' or 'petite' clearly differentiates the difference between noun and adjective or, in proceptual terms, between concept and process. The quantity is a *variable* and Cauchy says that 'cette variable *devient* ce qu'on nomme un infiniment petit' which translates to 'this variable *becomes* what one calls infinitely small.' In proceptual terms, therefore, Cauchy is speaking of a *process* that *becomes* small, rather than a *concept* that *is* small.

He continues by introducing +∞ as the limit of a sequence that increases indefinitely and –∞ as the limit of a sequence that decreases negatively. In this way his *quantities* now include the infinitely small and the infinitely large, although he continues to insist that *numbers* are always finite.

The numbers are extended to include positive and negative quantities in the form of signed numbers, but Cauchy also allows *variable quantities* (sequences of signed numbers) that may operate as infinitesimal or infinite quantities. It is a central point of contention among Cauchy scholars whether his continuum included only real numbers or whether it is enhanced by the inclusion of infinitesimals.

As a possible response to this dilemma, we note that if Cauchy sees points as lying *on* a line, in the manner of Greek geometry, then the line is not composed of points, but is rather an entity in itself on which points lie. Such points may be constants, staying in a fixed place, or variable quantities moving around. The question of the nature of Cauchy's continuum arises in the change in meaning, from Cauchy's view based on Greek geometry and the developments of Descartes, to Cantor's view of the continuum as the aggregate of its points, a viewpoint that is shared in our modern formulation of mathematical analysis and may cause a misinterpretation of Cauchy's original intentions.

At the opening of chapter II, Cauchy recalls his definition of infiniment petit in the following terms:

> On dit qu'une quantité variable devient infiniment petite, lorsque sa valeur numérique décroît indéfiniment de manière à converger vers la limite zéro.

His 'quantité variable infiniment petite' is a variable quantity given by a sequence of values whose absolute value (valeur numérique) converges to the limit zero; geometrically, this can be viewed as a sequence of points marked on a number line, successively closer to the origin.

### 4.3 Cauchy's conception of function

In chapter 1 of *Cours D'Analyse*, Cauchy introduces his notion of function, speaking explicitly about a relationship between two variables, say *x* and *y*, where one is the *independent variable x*, and the other is the *dependent variable y* whose values depend on *x*. He already has the familiar range of standard functions, including polynomials, rational functions, trigonometric functions, powers, exponentials and logarithms with a variety of behaviours that he needs to take into account. Throughout his text he speaks of functions in general terms, but his examples are always given as combinations of standard functions, which is to be expected in a textbook intended for student engineers.

He allows both explicit functions (where *y* is a function of *x*) and implicit functions (where *x* and *y* have a relationship between them). This enables him to consider a relationship such as $y = x^2$ where *y* is an explicit function of *x*, but when *x* is expressed in terms of *y*, there are two roots $y = \pm\sqrt{x}$ for positive *x*. He has his own special notion for such cases, using a single root $y = \sqrt{x}$ to denote the positive root $y = +\sqrt{x}$ and a double root sign for both roots, $y = \sqrt\!\!\!\sqrt{x}$. In the same manner, for the inverse sine function he uses $y = \arcsin(x)$ for the principal value in which $-\pi/2 \leq y \leq \pi/2$ and a special double parenthesis notation $y = \arcsin((x))$ to denote the multiple values. These conceptions should alert the modern reader to be aware that Cauchy uses terminology in a manner appropriate for his own era and this may differ significantly from modern usage.

**4.4 Cauchy's conception of limit**

Cauchy's conception of limit also allows it to have multiple values (Figure 2).

> Supposons, pour fixer les idées, qu'une variable positive ou négative représentée par $x$ converge vers la limite o, et désignons par A un nombre constant : il sera facile de s'assurer que chacune des expressions
>
> $$\lim A^x, \quad \lim \sin x$$
>
> a une valeur unique déterminée par l'équation
>
> $$\lim A^x = 1$$
>
> ou
>
> $$\lim \sin x = 0,$$
>
> tandis que l'expression
>
> $$\lim\left(\left(\frac{1}{x}\right)\right)$$
>
> admet deux valeurs, savoir, $+\infty$, $-\infty$, et
>
> $$\lim\left(\left(\sin\frac{1}{x}\right)\right)$$
>
> une infinité de valeurs comprises entre les limites $-1$ et $+1$.

Figure 2: Multiple values of limits

Here we see that the function $1/x$ has *two* possible limits as *x* tends to 0, which are $+\infty$, $-\infty$, and $\sin(1/x)$ has *an infinite number of limiting values* between $-1$ and $+1$. Cauchy's notation makes interesting reading for someone steeped in the concepts of modern analysis where a limit, if it exists, must be unique.

There is a subtle ambiguity in meaning of the opening sentence where he speaks of 'une variable positive ou négative représentée par $x$ converge vers la limite 0'. Does this mean that a single variable $x$ can move between being positive or negative (moving in from either side) or does it refer to two different cases: one negative moving in from the left and the other positive moving in from the right? In general he allows sequences to have terms with different signs, just as Leibniz did when proving the theorem that if $(a_n)$ is a monotone sequence of positive terms tending to zero then the alternating series $\sum (-1)^n a_n$ tends to a limit.

In dealing with the limit of $1/x$ where $x$ tends to zero, Cauchy separates out the two distinct cases where $\lim((1/x))$ has two values, $+\infty$ from the left and $-\infty$ from the right. This has possible subtle links to the meaning of the French term 'limite' which is used in a different way from the English term 'limit'. While both are used to mean the mathematical limit, the French term limite is also used to denote the endpoints of an interval. This brings in subtle shades of meaning in Cauchy's writing.

In the case of the two 'limites' for the variable $1/x$, we see that each one essentially adds a single endpoint to the curved lines making up the graph on either side of the origin. In exactly the same way, adding the limit points $x = -\infty$ to the negative part of the graph and $x = +\infty$ to the positive part of the graph complete the curved lines of the two parts of the graph by adding their 'limites'. Now both curved parts of the geometric graphs fit the Greek notion of a line with endpoints at each end.

The final case he mentions concerns the limit of $\sin(1/x)$ having 'an infinity of values between the limits 1 and +1'. This proves to be even more interesting and we will return to it at the end of the next section when we have discussed the notion of continuity.

### 4.5 Cauchy's conception of continuity

Cauchy presents his definition of continuity in three consecutive stages, one following immediately after another. (Figure 3.)

§ II. — *De la continuité des fonctions.*

Parmi les objets qui se rattachent à la considération des infiniment petits, on doit placer les notions relatives à la continuité ou à la discontinuité des fonctions. Examinons d'abord sous ce point de vue les fonctions d'une seule variable.

Definition 1 {
Soit $f(x)$ une fonction de la variable $x$, et supposons que, pour chaque valeur de $x$ intermédiaire entre deux limites données, cette fonction admette constamment une valeur unique et finie. Si, en partant d'une valeur de $x$ comprise entre ces limites, on attribue à la variable $x$ un accroissement infiniment petit $\alpha$, la fonction elle-même recevra pour accroissement la différence

$$f(x+\alpha) - f(x),$$

qui dépendra en même temps de la nouvelle variable $\alpha$ et de la valeur de $x$. Cela posé, la fonction $f(x)$ sera, entre les deux limites assignées à la variable $x$, fonction *continue* de cette variable, si, pour chaque valeur de $x$ intermédiaire entre ces limites, la valeur numérique de la différence

$$f(x+\alpha) - f(x)$$

décroît indéfiniment avec celle de $\alpha$. } Definition 2 { En d'autres termes, *la fonction $f(x)$ restera continue par rapport à $x$ entre les limites données, si, entre ces limites, un accroissement infiniment petit de la variable produit toujours un accroissement infiniment petit de la fonction elle-même.*

Definition 3 {
On dit encore que la fonction $f(x)$ est, dans le voisinage d'une valeur particulière attribuée à la variable $x$, fonction continue de cette variable, toutes les fois qu'elle est continue entre deux limites de $x$, même très rapprochées, qui renferment la valeur dont il s'agit.

Enfin, lorsqu'une fonction $f(x)$ cesse d'être continue dans le voisinage d'une valeur particulière de la variable $x$, on dit qu'elle devient alors *discontinue* et qu'il y a pour cette valeur particulière *solution de continuité*.
}

Figure 3: Cauchy's definition(s) of continuity

His first definition speaks of a function $f(x)$ which, for each value of *x* between two given limits, has a unique finite value. The term 'limite' here means the endpoints of an interval. A function is said to be continuous if, for a value *x* between these limits, an infinitesimal increase α in *x* causes the difference $f(x+\alpha)-f(x)$, which depends on the new variable α, to decrease indefinitely with α. Note that this a *process* definition: the difference *decreases indefinitely* as α decreases.

The second definition is given in italics, which suggests that this is his main definition that he may wish his students to commit to memory. It states that (for *x* between the given limits), an infinitesimal change in *x* gives an infinitesimal change in the function $f(x)$. Notice the subtle change in meaning as the *process* of change, formulated in the first definition 'as the difference decreases indefinitely', subtly becomes a *concept* of change using the noun phrase 'infiniment petit' in the second.

Of course, Cauchy himself is generally careful to explain that a null sequence *becomes* infinitely small, rather than saying that it *is* an infinitely small quantity. However, this façon de parler introduces not one, but two, new elements in the definition of continuity.

First it gives an operable definition of computing the difference $f(x+\alpha) - f(x)$ for an infinitesimal α, because α is generated by a sequence $a_1$, $a_2$, $a_3$, …, $a_n$… of finite quantities (numbers with signs); see Section 6 above. Thus, the sequence of values $f(x+a_1)$, $f(x+a_2)$, $f(x+a_3)$, … is a computable sequence of values which can be checked to determine if it tends to zero. (We note that nowhere in the book does Cauchy perform a specific numerical calculation. However, such a computation would often become quite complicated and we interpret Cauchy's presentation as taking the pragmatic route of giving specific examples where the general idea of convergence of the terms appears to be evident without the need to do what might be a tricky calculation.)

Secondly, and more significantly, the language compresses the process of becoming arbitrarily small from the phrase 'decroit indefiniment' into the *noun* 'infiniment petit', as in Cauchy's 1823 definition cited in Section 6. Such a noun can now function as a mental entity as if it were a manipulable concept in its own right.

Cauchy does not remark explicitly on this transition from process to mental entity, however, it is a natural development, as expressed succinctly by Laugwitz:

> Having recognized the usefulness of infinitely small magnitudes in research, he must have been tempted not only to use them in concurrently written textbooks but also to justify them rigorously. This was a stepwise process. At the beginning we still encounter the traditional locutions about variables with limit zero or infinity. In time, the infinitely small magnitudes, soon referred to as numbers, acquire independence and are handled like 'genuine' numbers. (Laugwitz 1997, p. 657.)

Such a manner of speaking is common in the informal conversation of mathematicians, and also in the increasingly subtle thinking of growing children as they learn to switch flexibly from thinking of an operation as a process to be carried out and the value given by the process (Gray & Tall, 1994, Tall, 1991).

The third definition moves on to expand the definition of continuity to include the notion of discontinuity. This says that a function is continuous *in the neighbourhood of a particular value* of the variable *x*, whenever it is *continuous between two limits for x, even very close* (our italics). Cauchy contrasts this with the notion of a function being discontinuous at a particular value of *x* if it ceases to be continuous in a neighbourhood of *x*. In particular, he notes that, in a neighbourhood of zero the function $1/x$ 'becomes infinite', and, as a consequence is *discontinuous* there.

These definitions of continuity may be read in various ways depending on the experiences of the commentator. For example, they may be read as introducing the formal definition of pointwise continuity. Yet even though Cauchy speaks of continuity *at a particular value* of *x*, he interprets this *in the neighbourhood* of that value *between two limits for x, even very close*. This is consistent with a level two conception of continuity, visualizing a perceptually continuous function being drawn between two limits. It is consistent with Cauchy's idea of a *variable* quantity moving along the curve. Hence it relates to the continuity of (a part of) the curve as a variable point moves along it, rather than as a pointwise definition of continuity found in modern analysis.

The possible meanings of the definitions of continuity and limit become clearer when we consider the examples that Cauchy gives in *Cours D'Analyse* (figure 4).

In this figure, the functions in the first list are said to be continuous 'in the neighbourhood of a finite value (attributed to) $x$' where that finite value lies between the 'limites' $-\infty$ and $+\infty$. A modern mathematician may translate this as saying 'for all finite $x$'. However, Cauchy sees these points varying *between* the 'limites', in a manner which is quite natural when those limits involve potentially infinite quantities either horizontally or vertically at either end.

The function $a/x$ is continuous between two separate 'limites', one between $-\infty$ and $0$, the other between $0$ and $\infty$.

Par suite, chacune de ces fonctions sera continue dans le voisinage d'une valeur finie attribuée à la variable $x$, si cette valeur finie se trouve comprise :

Pour les fonctions
$a + x$
$a - x$
$ax$
$A^x$
$\sin x$
$\cos x$
entre les limites $x = -\infty$, $x = +\infty$;

Pour la fonction
$\dfrac{a}{x}$
1° entre les limites $x = -\infty$, $x = 0$,
2° entre les limites $x = 0$, $x = \infty$;

Figure 4: Continuous functions

Returning now to the limit of $\sin(1/x)$ as $x$ tends to zero considered earlier in figure 2, we find a much more interesting phenomenon. Once again the function is defined and continuous in a neighbourhood of each finite value of $x$ to the left and to the right of the origin. Following the example of $a/x$, we can again see the function is continuous in two parts, '1st between the limits $x = -\infty$ and $x = 0$' and '2nd between the limits $x = 0$ and $x = \infty$'. The limit $x = -\infty$ may be adjoined to the left portion of the graph and $y = \infty$ to the right, consistent both with the Euclidean notion of a point being an endpoint of a (curved) line and also the French terminology of 'limite' as the endpoint of an interval.

This follows the same pattern as before, adjoining a single extra point at the end of a line. However, near the origin, a different phenomenon occurs. (Figure 5.)

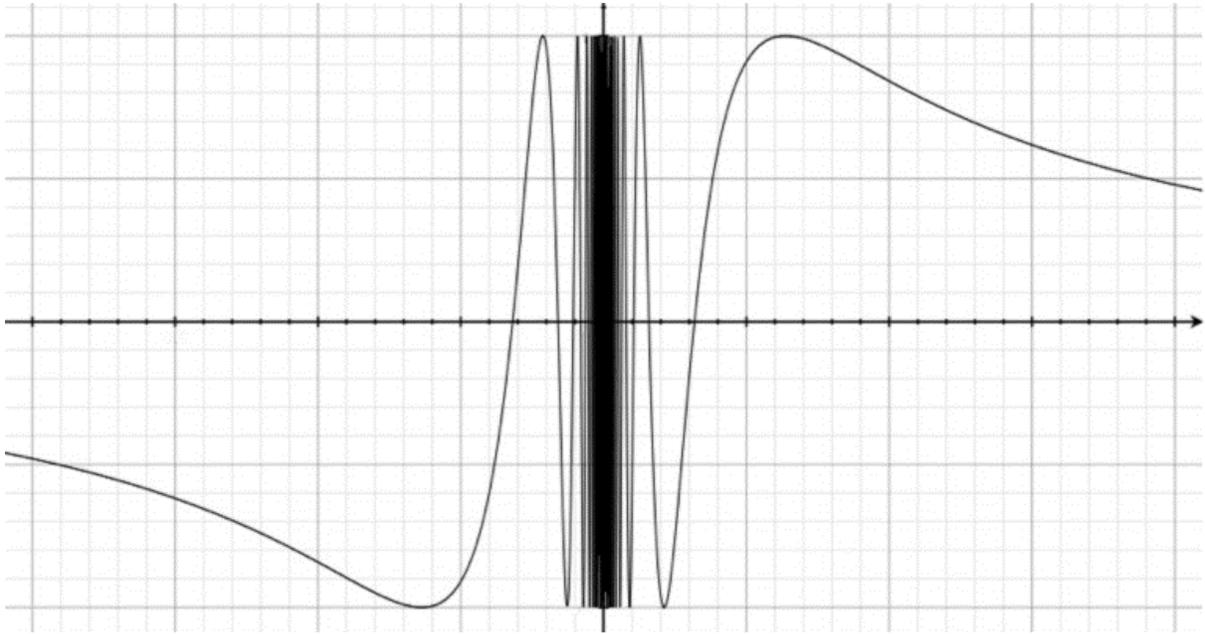
Figure 5: The graph of sin(1/*x*)

If we restrict ourselves either to moving in from the left or from the right, drawing the curve with a physical pencil, our strokes move closer and closer to the vertical interval on the *y*-axis from –1 to +1 and we end up running our practical pencil up and down closer and closer to the *y*-axis which acts as a barrier, as a limit to our progress.

Figure 6 show the curve being drawn from the right as *x* moves down to zero. It will eventually end up moving up and down perceptually between –1 and 1 on the *y*-axis, even though, theoretically it is close to it rather than being on top of it. At Donald's level two of perceptual continuity, we sense the hand moving up and down, covering the *y*-axis within the thickness of the pencil on the interval from –1 to +1.

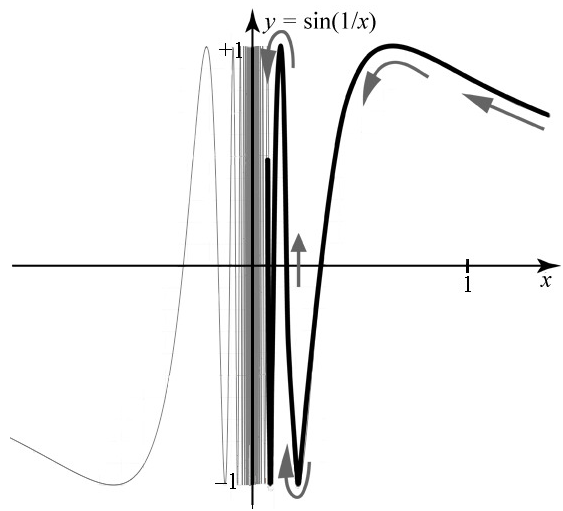
Figure 6: A smooth pencil drawing moving from the right to the left

The picture becomes more interesting when we follow Cauchy's idea of a focusing on a sequence $a_1, a_2, a_3, \ldots, a_n, \ldots$ where the x-values tend to zero. Figure 7a, shows the first twenty terms when the nth term is calculated for $x = 1/n$. The terms successively move right to left along the curve. Figure 7b, shows the terms from $n = 100$ to $1000$. The finite width of the dots marking the points give an infinite number of points (apparently) on the vertical axis from –1 to +1.

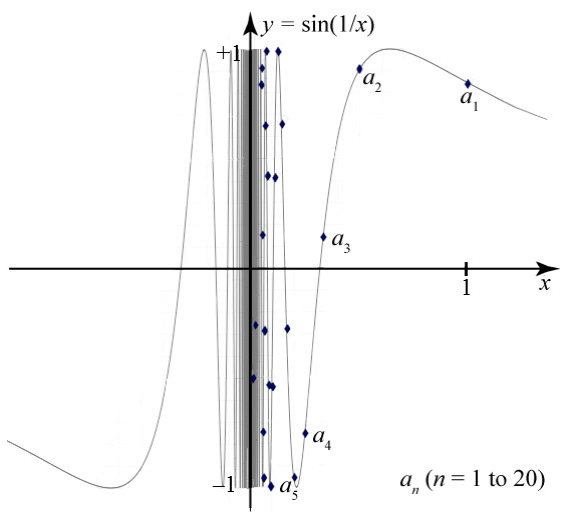 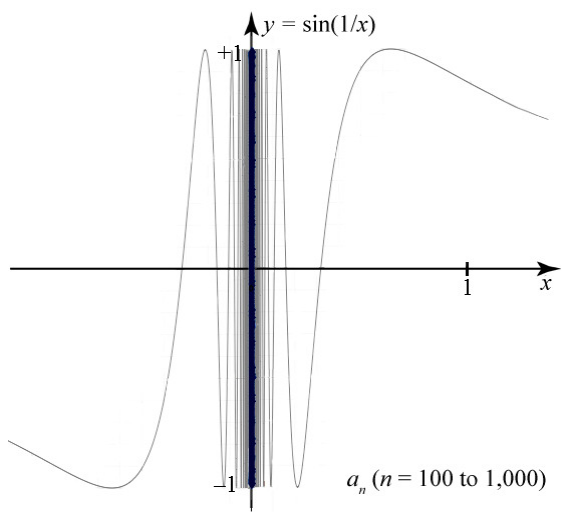

Figure 7a: $a_n$ ($n = 1$ to 20)    Figure 7b: $a_n$ ($n = 100$ to 1,000)

We now *see* the limiting values, in a sense which is compatible with Cauchy's idea that the limit of $\sin(1/x)$ as x tends to 0 is, in Cauchy's own words:

  une infinité de valeurs comprises entre les limites –1 et +1.

### 4.6 Cauchy's development of meaning 1821–1853

In his *Cours D'Analyse*, Cauchy identifies quantities as (signed) lengths. An 'infiniment petit' is not a fixed length, it is a variable quantity where a succession of values *becomes* arbitrarily small. It occurs in expressions such as $f(x+\alpha)$ which represents a succession of values as α becomes small.

However, he uses the term 'infiniment petit' *as a noun* and computes with it as a variable quantity that 'becomes' small. In this way it takes on a cognitive existence of its own. He later uses infinitesimals not only in computing limits, derivatives and integrals, but as explicit entities in their own right. For instance, in 1827, he uses quantities involving infinitesimals as the limits (endpoints) of integrals in more complicated expressions such as

$$\sum_{n=1}^{n=\infty} \int_{\frac{2n\pi}{\alpha}-\delta}^{\frac{2n\pi}{\alpha}+\delta} \frac{\varepsilon\,\varphi(r)\,dr}{\varepsilon^2 + (1-\varepsilon)\left(2\sin\frac{\alpha r}{2}\right)^2}$$

where α, ε, δ are infinitesimals and ε is an infinitesimal of the second order compared to δ (for instance $\varepsilon = \delta^2$), (Cauchy 1827, p. 188). He also uses these highly complex expressions in his development of a formula for a unit-impulse, infinitely tall, infinitely narrow delta function defined in terms of infinitesimals. In this way, he invented the Dirac delta function in infinitesimal terms a century before Dirac. (See Laugwitz, 1989, p. 230.)

In 1853, he responded to criticisms of his theorem, which states that 'the limit $s(x)$ of a convergent sequence of continuous functions $s_n(x)$ is continuous', and is now seen to contradict ideas relating to Fourier series. Cauchy is fully aware that the convergence condition fails for the Fourier series

$$s_n(x) = \sum_{k=1}^{n} \frac{\sin kx}{k}$$

if the hypothesis of convergence is interpreted to apply only at ordinary (real) points. Cauchy responds by reasserting his original definitions of continuity and convergence, calculating the error term $s(x) - s_n(x)$ for the variable quantity $x = 1/n$, and concluding that the convergence hypothesis is *not* satisfied. What Cauchy is pointing out is that, if the hypothesis of convergence is interpreted as applying at additional values (in particular, infinitesimal ones generated by $x = 1/n$), then the sum theorem remains valid (Cauchy, 1853, p. 33). See (Bråting, 2007) for details.

This fine distinction has led to enormous controversy as later authors attempt to interpret Cauchy's ideas through a conflation of his original worded definition and subsequent set-theoretic definitions of pointwise and uniform convergence. In a similar vein, Núñez et al. (1999) bring together the names of Cauchy and Weierstrass to speak of the 'Cauchy-Weierstrass definition of continuity'. J. Gray (2008, p. 62) lists continuity among concepts Cauchy allegedly defined using 'limiting arguments', but as we discussed in Section 9, 'limits' appear in Cauchy's definition only in the sense of endpoints of an interval, rather than 'limits' as in variables tending to a quantity. Not to be outdone, Kline (1980, p. 273) claims that "Cauchy's work not only banished [infinitesimals] but disposed of any need for them." Hawking (2007, p. 639) does reproduce Cauchy's infinitesimal definition, yet on the same page claims that Cauchy "was particularly concerned to banish infinitesimals," apparently unaware of a comical non-sequitur he committed.

Borovik & Katz (2012) argue that the term 'Cauchy-Weierstrass definition of continuity' is inappropriate, because Cauchy's definition is based on infinitesimals rather than on Weierstrassian epsilon-delta methods. They note that Cauchy's notion of limit was kinetic rather than epsilontic (see also Sinaceur 1973), and that his work on the orders of growth of infinitesimals inspired later theories of infinitesimals by Paul du Bois-Reymond, as noted also by Borel (1902, p. 35–36).

Extending Bråting's analysis of Cauchy's infinitesimals, Katz & Katz (2011) argue that the traditional approach to Cauchy following Boyer and Grabiner contains an internal contradiction: on the one hand, Cauchy's infinitesimal definition is interpreted as applying only at standard real values of the input; yet the condition of Cauchy's sum theorem of 1853 is routinely interpreted as a condition of uniform convergence, which is only possible if one allows variable (in particular, infinitesimal) values of the input (such as $x = 1/n$).

Katz & Katz (2012) argue that post-Weistrassian infinitesimal-free readings of Cauchy focus on those aspects that fit the modern use of epsilon arguments without fully taking into account his use of infinitesimals. Błaszczyk et al. (2013) re-examine Grabiner's claim that Cauchy's proof of the intermediate value theorem contains germs of the epsilon-delta method and offer an alternative view that Cauchy's proof fits squarely within his infinitesimal approach.

He is using the same modus operandi from his earlier years, computing functions whose inputs are variable quantities, in this case the 'infiniment petit' variable $(1/n)$ which is specified as a sequence of ordinary numbers with limit zero.

He does not see the number line as the aggregate of real numbers as in the modern set-theoretic viewpoint. Rather, he sees a line with points *on* it. He imagines quantities on this line where a quantity may be a fixed constant, or it may vary as in the sequence $(1, ½, …, 1/n, …)$ whose value is computed at each fixed number in the sequence to give a new sequence $(f(x+1), f(x+½), … f(x+1/n), …)$ of numbers which he sees once more as a variable quantity.

Cauchy's work can be seen more equitably as a natural development of the evolution of mathematical ideas in a Lakatosian sense (Koetsier, 2010). Cauchy makes sense of infinitesimals in a concrete fashion in terms of sequences of signed numbers that tend to zero. He then deals with increasingly sophisticated situations in which his infinitesimals as processes of becoming small are mentally compressed into infinitesimals as manipulable concepts.

It should not be assumed that the evolution of mathematical ideas always produces 'better' results or that our current views are final and will not continue to evolve. For instance, our current set-theoretic ideas of functions and relations lose some subtleties that exist in Cauchy's version of calculus. In particular, the pragmatic decision to make functions single-valued loses the dynamic idea of two variables 'co-varying' which is part of the current development of educational views of the calculus. If we consider, for example, the inverse function of $\sin x$, we are now forced to artificially choose a part of the domain, for instance selecting $\sin^{-1}(x)$ to lie between $-\pi/2$ and $\pi/2$. But Cauchy imagined a multi-valued function $\arcsin((x))$ and, if we were to choose any one of these values and denote it by $y$ where $y = \sin x$ and allow it to vary on the graph of $\sin x$, then Cauchy's notion of continuity would allow this particular value of *y* to co-vary with *x* as *x* moves over an artificially selected endpoint such as $x = \pi/2$.

More serious difficulties arise with relations such as $x^2 + y^2 = k$, which is a circle and *y* is no longer considered to be a function of *x*. Yet, if one traces a finger round such a circle, the values of *x* and *y* co-vary in the sense of Cauchy. This arises quite naturally in physical situations such as the motion of a weight hanging on a string from a fixed point and moving in a horizontal circle where the Cartesian coordinates of its position in a horizontal plane co-vary. The modern single-valued notion of function has many advantages but it does not fit all possible cases, particularly in naturally occurring 'continuous' change as modelled in applied mathematics. The ideas of Cauchy continue to have merit, even if we insist on new set-theoretic definitions of functions and continuous change in mathematical analysis.

## 5. THE DEVELOPMENT OF FORMAL PROOF

The shift to the formal world of Hilbertian set theory involves a fundamental change of meaning in which the concepts are not naturally occurring phenomena that 'have' properties, as in Euclidean geometry or in arithmetic and its generalization in algebra. In today's formal presentation they are set-theoretic concepts given as a list of axioms and what matters are the consequences of the properties specified explicitly by the axioms and any additional definitions within the axiomatic system.

An axiomatically defined structure may apply in many quite different situations. Often, as in the case of the axiomatic definition of a group there are many different examples of groups that can be partitioned into collections that are isomorphic to each other. In the case of a 'complete ordered field' there is, up to isomorphism, only one example: the real numbers, which may be embodied as a complete number line or expressed symbolically as infinite decimals with appropriate properties.

As Cantor extended numbers both in terms of the complete ordered field of real numbers and his theory of infinite cardinals, he believed that he had not only given a logical foundation to real analysis, he had simultaneously eliminated infinitesimals, declaring them to be the *cholera bacillus* of mathematics (in a letter of 12 December 1893, quoted in Meschkowski 1965, p. 505). Successive generations of pure mathematicians have followed his lead and accepted epsilon-delta analysis as *the* proper formal approach.

But the evolution of mathematical ideas moves on. The formalism of Hilbert does not grant the right to Cantor or anyone else to deny the existence of infinitesimals on logical grounds.

Cantor's rejection of infinitesimals on the real number line relates to his perception of real numbers, which form a complete ordered field and therefore, rightly, cannot include infinitesimals. However, this does not mean that there cannot be ordered fields *K* that *extend* $\mathbb{R}$ that contain not only real numbers, but also positive elements *x* that satisfy $0 < x < r$ for all positive $r \in \mathbb{R}$. Such an element is a positive infinitesimal in the larger ordered field *K*. There are many examples of such fields, including not only the hyperreal numbers in non-standard analysis, but also simpler fields such as the field of rational functions in an element *x* where we define $0 < x < r$ for any positive real number *r*, or the field of quotients of power series in *x*, which Tall (1979) showed was sufficient to deal with the calculus of all analytic functions that occur in regular calculus courses, including the calculus of Leibniz, without the need for the sophisticated logic of non-standard analysis.

If *K* is any ordered field that contains the real numbers as a proper subfield, then a simple structure theorem using the completeness of $\mathbb{R}$ proves formally that every finite element of *K* (lying between two real numbers) is either a real number or a unique real number plus an infinitesimal (Tall, 1982b, 2001). Although a real number *c* cannot be distinguished by the naked eye from a nearby element $c+\varepsilon$, when $\varepsilon$ is an infinitesimal, a linear map of the form $f(x) = mx + b$ can be used to magnify infinitesimal detail. For instance, the map $\mu: K \to K$ given by $\mu(x) = (x-c)/\varepsilon$ maps *c* to 0 and $c+\varepsilon$ to 1, with any other nearby point $c+k\varepsilon$ mapped on to *k* (Stroyan, 1972; Tall, 1980). In particular, for any real number *k*, this allows infinitesimal detail $c+k\varepsilon$ to be mapped on to the whole real number line. Of course, this map will transform the point $c+\varepsilon^2$ onto $\mu(c+\varepsilon^2) = \varepsilon$ and the difference between the images is again too small to see with the perceptual eye, because higher order differences remain infinitesimal after the act of magnification (see Tall (2001, 2013).

This simple theorem has a remarkable interpretation. If we call the elements of *K* 'quantities' and the elements of $\mathbb{R}$ 'constants', then the ordered field *K* contains infinite quantities, whose inverses are infinitesimal quantities, and the finite elements in *K* are either constants or a constant plus an infinitesimal. In this way *any* ordered extension of the real numbers defined formally can be verbalised in the language of Cauchy! While this may not in any way 'explain' Cauchy's thinking, it shows that his infinitesimal ideas can evolve to take their rightful place in today's axiomatic formal approach.

Thus infinitesimals, so strongly denied by Cantor, are a natural consequence of Hilbertian formalism. While we should still honour Cantor for putting mathematical analysis on a rigorous basis by defining the real numbers as a complete ordered field, we now realise that he did not eliminate infinitesimals at all. They simply lie in a larger mathematical structure that can be defined formally, manipulated symbolically, and visualised perceptually.

It therefore does no service to argue whether Cauchy was a proto-Weierstrass foreseeing the development of epsilon-delta analysis or a proto-Robinson foreseeing the advent of non-standard analysis. He was neither. He was an intellectual practitioner who lived in an era building on Greek ideas of geometry extended through Descartes' links between geometry and algebra to operationalize the notion of continuity, with a cultural vision of 'limites' and infinitesimals that were generated by ordinary sequences of numbers with limit zero.

Rather than analysing his techniques from a modern viewpoint, it would be more appropriate to seek to understand Cauchy's ideas as part of the evolution of mathematical ideas in the style of Lakatos (1976). We will then recognize his unique contribution to mathematics based on a combination of human perception, operational symbolism and intellectual reason using ways of thinking that are appropriate to his own era.

**5.1 The evolution of theoretical frameworks**

Since Cauchy's time, new ways of thinking have emerged that subtly modify the ways in which we think about mathematics. Cantor's view of the number line emerged at a time when axiomatic definitions and formal proof were evolving, yet still implicitly attached to human perception of the number line and the arithmetic of numbers. The full implications of Hilbert's vision of formal axiomatics had yet to be fully grasped.

Hilbert's views were to flower into the set-theoretic formulation of definitions and proof found in modern mathematical analysis. Yet is was not long in development before Gödel's theorem showed that formal systems that included the integers would have too many theorems to be able to prove each of them in only a finite number of steps. Then came the logical ideas of Robinson that placed the real number system within a larger hyperreal number system with infinitesimals that again led to vigorous controversy between different communities of practice.

As we consider these ongoing disputes, we recall the scheme of Perry (1970), who described the development of college students setting out with a belief that truth exists in absolute terms of right and wrong, so that some authorities are correct and others are not, passing through a stage that realized there may be alternative answers where some are better than others, and moving to the recognition that there are different, legitimate alternatives.

As we review the varying opinions of Cauchy's ideas it becomes evident to us that it is naïve to see differing viewpoints simply as right or wrong in absolute terms and to acknowledge that various principles may operate more successfully in one context as compared with another. In particular, in attempting to make sense of the intellectual developments of students, we would be wise to acknowledge explicitly the evolving nature of our own theoretical frameworks, in particular, the shift from natural thinking about quantities that we imagine as being 'arbitrarily small' and the formal epsilon-delta definitions of mathematical analysis. This would enable us to have a broader view of the conceptions of others that may make complete sense in one community of practice yet be less appropriate in another. It continues to encourage pure mathematicians to pursue the rigour of formal epsilon delta analysis while allowing applied mathematicians to utilise the powerful blend of dynamic human perceptions and the symbolic use of operational symbolism to formulate and predict solutions of problems in applications.

It may even allow us to move to the final level of Perry's scheme where we allow the possibility of different, equally legitimate alternatives, where each applies appropriately in a particular theoretical or practical situation.

## 6. STUDENT CONCEPTIONS

Cauchy's great insights into the nature of continuity, limits and infinitesimals reveal a refined human mind seeking to make sense of ideas of variability in operational terms that sought to make sense of the potentially infinite process of 'getting close' to a 'limiting value'. Our students today may not have Cauchy's exceptional skills, but they do share with him a human brain that perceives fundamental concepts of variation and continuity in the same biological manner.

Roh (2008) introduced an innovative approach to the formal definition of limit of a sequence $a_1$, $a_2$, ... using human perception by placing a transparent strip with parallel lines distance $\pm\varepsilon$ apart over a picture of a sequence of points plotted as a graph to see if the sequence lay between the lines for all but a finite number of values. By imagining this activity for smaller values of ε, the students were encouraged to suggest their own definition of convergence. The research revealed three distinct conceptions of a limit:

(i) an *asymptotic* view in which the sequence approaches but does not reach the limit, such as (1/*n*) 'approaching' zero, but excluding a constant sequence which do not approach the limit because it is already there;

(ii) a *cluster point* view in which the sequence may cluster round several values, as happens with $(-1)^n(1+1/n)$);

(iii) *the modern limit concept* with a unique limit.

Remarkably, all these three conceptions arose in history. The first—which often occurs in our students—is found in the original conception of Leibniz and also in Newton's notion of *prime* and *ultimate* ratios in which the prime ratio of $(x+o)^n - x^n$ to *o* for a non-zero value *o* tends to the ratio $nx^{n-1}$ to 1, but is never equal to it (Newton, 1704). The second relates to the notion of limit expressed by Cauchy in the *Cours D'Analyse*, where sequences cluster around one or more points that may be infinite in size or infinite in number. The third is the modern definition that replaces a sense of dynamic change by a logical challenge.

This analysis shows clearly that what are often described in the research literature as student 'misconceptions' are more appropriately described as 'pre-conceptions' that occur at early stages of development and may or may not blossom into a more widely accepted formal concepts.

Like Cauchy, our modern students have experiences with drawing graphs and manipulating symbols but are unlikely to have encountered the formal set-theoretic epsilon-delta approach. We are confronted with the serious question: why, if the notion of dynamic change worked for a mathematician as good as Cauchy and continue to work in applications today, do we introduce all our students in their early experiences in the calculus with a definitional version of limit that is so different from their current experience?

The mathematical community shares a common expertise acknowledged by mathematics educators in which the limit concept has allowed them to produce far more powerful methods of definition and proof and it is natural for them as experts to wish to introduce their novice students to this more powerful mode of operation. Some students do indeed benefit by a 'formal' approach based on 'extracting meaning' from definitions and deductions rather than a more 'natural' approach, 'giving meaning' to definitions based on their previous experience of embodiment and symbolism. Pinto & Tall (1999, 2001) used such a framework to analyse the progress of students learning mathematical analysis in a standard university course and found that students could be successful both with either a 'natural' or 'formal' approach, but that many found conflict between the new ideas and their previous experience that often led them into difficulties. Weber (2004) expanded this observation by noting that many students followed 'natural' or 'formal' routes but many more solved their difficulties by a 'procedural' route in which they side-tracked their problems by learning proofs by rote to reproduce in examinations.

Alcock and Simpson (2005) reported the successful learning of students in an analysis course who operate 'almost exclusively by means of verbal and algebraic reasoning, and tend not to incorporate visual images into their work'. This is consonant with acknowledging a spectrum of approaches using natural and formal ideas. Indeed, it is an advantage in mathematical research for individuals to have different ways of thinking about problems to share different aspects of the creation of new theoretical developments.

Mejia-Ramos and Tall (2004) put forward the thesis that calculus naturally belongs to the worlds of embodiment to conceptualise problem situations and symbolism to formulate mathematical models using the operations of calculus to solve problems. This is the foundation that is required in real world problem-solving in applications. According to this view, the calculus belongs to the natural world of human perception and operation, using familiar pictures and symbolic calculation and manipulation, while analysis belongs to the formal world of definition and proof. It is natural to encourage those who use the calculus in applications to be aware of its natural strengths and possible weaknesses. Meanwhile, a formal approach can involve standard analysis with the Weierstrassian definitions, or non-standard analysis using the logical notions of infinitesimals, both of which may be seen to evolve through succeeding generations using a subtler formal interpretation of the conceptions of limit and infinitesimals envisaged by Cauchy.

These observations relating to student development are consonant with the historical evidence of conceptual changes that occurred over the years, from the separation between the platonic geometry of the Greeks and the arithmetic of number, through the evolution of operational symbolism of various kinds of number, from magnitudes without sign, to signed numbers, rational and irrational numbers, variable quantities that may become large or small, and other developments in complex numbers and beyond. There is also the fundamental schism between mathematics that we sense in the external world and imagine as perfect platonic ideas and the subtle shift to axiomatic formal mathematics based on set-theoretic definitions.

These are reflected in the theoretical framework built from the cognitive psychology of human perception and dynamic action in Donald's level 1 and 2, moving on to level 3 extended awareness as mathematics develops through increasingly sophisticated conceptual embodiment and operational symbolism and later makes the significant transition into axiomatic formalism.

## 7. IMPLICATIONS FOR EDUCATIONAL THEORISTS AND MATHEMATICS EDUCATORS

This analysis of the mathematics of Cauchy seen through the ideas of those who have the benefit of the later evolution of ideas has aspects in common with the analyses of student thinking by experts in mathematics and mathematics education who have more sophisticated views that they wish to convey to students in the calculus and later in mathematical analysis.

There are many possible approaches including the standard 'intuitive limit' approach in modern College Calculus texts such as Stewart (2011), Larson & Edwards (2009), a multi-representational approached based on numeric, symbolic and graphic representations (Hughes-Hallett *et al.*, 2009), the cooperative learning approach to quantified definitions in Dubinsky's APOS theory (Cottrill et al. 1996), or an infinitesimal approach based on non-standard analysis by Keisler (1976). Computer technology introduces completely new possibilities, building on human perception and action, passing much of the technical operation to software such as Mathematica, Maple, Cabri Géomètre, SketchPad, Geogebra, on computers, graphic calculators, and the enactive interface of iPads, interacting with the internet and learning cooperatively in classrooms, using wireless technologies and other modes of representation and communication.

Faced with such a plethora of possibilities, how does one cut through the morass of detail and seek the fundamental ideas in the calculus that both make sense to students and also provide them with the computational power required in today's modern technological world? The analysis of both historical and cognitive evolution of mathematical concepts given here suggests the need to reflect on the manner in which human thinking grows from perception and action into imaginative thought experiments and operational symbolism, developing increasingly sophisticated forms of human reason.

This confirms the clear schism between the dynamic calculus of Cauchy and the set-theoretic formal approach of modern analysis as revealed by the thorny path of the evolution of mathematical concepts in history and the difficulties encountered by students today.

On the one hand the approach of Dubinsky and his colleagues (Asiala et al., 1996) sought to encourage students to make sense of quantified statements through cooperative programming. However, even though many students were able to conceptualise functions as processes, only a few understood the formal definition of limit, and not one student in a study by Cottrill et al. (1996) constructed the formal definition spontaneously.

On the other hand, Núñez et al. (1999) saw 'formal epsilon-delta continuity' as being very different from 'natural continuity' based on human embodiment and perception. Lakoff & Núñez (2002) took the argument further, rejecting the views of mathematicians that 'mathematics is universal, absolute and certain' and contrasted the 'beautiful story of mathematics' presented by mathematicians with 'the sad consequences' that the formal approach 'intimidates people' and 'makes mathematics beyond the reach of even intelligent students with other primary interests and skills,' (pp. 339–341).

Yet some students, perhaps only a few, espouse a formal approach (Alcock & Simpson, 2005) while many others find the transition from natural experiences of visual and symbolic school mathematics to the formal limit to be problematic.

Instead of seeing one approach to be 'good' and the other 'bad', depending on one's point of view, it is more appropriate to construct a mature understanding that acknowledges there are alternative approaches that are valid in their appropriate contexts and that some approaches may be more appropriate than others in a given context.

With the premise that our conceptions depend on our experience, it may help to reflect on the idea that what is 'natural' may depend on how each of us develops over time. The calculus of Cauchy imagines natural variation of quantities that have arisen through two millennia of observing straight lines and smooth curves in Greek geometry and the development of calculus where a smooth curve is a mental extrapolation of the movement of a finger drawing a curve using a pencil on paper or a stick drawing a line in the sand. Initially, when these experiences were symbolised, it was always in terms of regular functions given by algebraic and trigonometric expression. Now we have a broader concept of functions which may be both continuous yet so irregular that they are nowhere differentiable.

To take account of the human origins of mathematical thinking, Tall (2012) proposed an approach based on perceptual change that sees continuity in terms of the dynamic continuous movement of a point (at Donald's level 2) and differentiability in terms of 'local straightness' in which a curve, when highly magnified, looks increasingly like a straight line. This takes us back to the world of Leibniz who imagined his curves as polygons with an infinite number of infinitesimal straight sides.

A 'locally straight' view provides an aesthetic extension to Leibniz's vision, in that his infinitesimally sided polygons have corners. These may only involve an infinitesimal change in direction, yet they are corners all the same. The 'locally straight' approach looks closely at a smooth graph and sees it looking locally like a straight line *everywhere*. This links to Cauchy's vision of a graph as a curved line on which a variable point may vary in a dynamically continuous fashion; local straightness also reflects the 'natural' idea of continuity conceived by Lakoff and Núñez in which the graph has a tangent at every point.

A locally straight approach also allows the individual to imagine the difference between a function simply being continuous and one that is suitably smooth for the purposes of calculus. A non-differentiable function no longer magnifies to look straight everywhere. It may have a corner at a single point, or it may be so wrinkled that it *nowhere* magnifies to look straight (Tall, 1982a).

This radically challenges the view of what is 'natural' according to Lakoff and his colleagues. If one can imagine a function that is very wrinkled, then being non-differentiable becomes a 'natural' concept. It is all a matter of building on appropriate experience.

The denigration of formal approaches of mathematicians and its 'sad consequences' described by Lakoff and Núñez remains true if one only considers traditional approaches to learning mathematics, but it should now be considered in a broader light. Even though many students find difficulty with the formalities, those who learn to work with formalism are also human beings using human brains. Formal mathematics develops a network of formal theorems, some of which prove to be *structure theorems* that reveal formal structures possessing embodied and symbolic meanings. (Tall, 2001, 2013). Typical examples are:

> A finite group specified formally by axioms is isomorphic to a group of transformations that has an embodied representation as the symmetries of a figure and a symbolic representation as a subgroup of the permutation group $S_n$;

> a finite *n*-dimensional vector space over a field *K* is isomorphic to $K^n$, which provides a symbolic interpretation using matrices and, in the case where $K = \mathbb{R}$ and *n* = 2 or 3, it also provides an embodiment in two or three dimensional space;

> a formal complete ordered field has a unique embodiment as a geometric number line with points representing rationals and irrationals that may also be represented symbolically as the arithmetic of decimal numbers.

With the possibility that formal axiomatic systems lead back to embodied and symbolic representations now enhanced by formal proof, we see that 'the beautiful story of mathematics' is indeed an aesthetic human creation in which formal ideas may also be imagined and manipulated by combining human perception, operational symbolism and axiomatic formal proof. What happens in practice is that mathematicians develop personal preferences for one or more of these aspects.

The moral of this story is that our cultural views of mathematics, even if they are widely shared, such as the principle of building the calculus on an 'intuitive' version of the limit concept, are fundamentally dependent on our shared history of personal experiences. While we persist in teaching our students ideas about limits that are implicitly founded on an axiomatic formal approach to mathematics, we produce a succession of generations that pass on ideas that may appeal to a few but are confusing for many. This is often the case in the teaching of college calculus where students are required to work through compendious volumes that cover every possible topic.

There is a need to cater for a wider spectrum of students who furnish a broader range of participation in our complex society. We certainly require mathematicians who can think formally and evolve new theories of mathematics who benefit from a formal approach. We also require pragmatic mathematicians who apply mathematical ideas to solve real world problems. And we need to provide courses in the calculus that make genuine sense of ideas of the dynamic rates of change and growth of variable quantities for a much broader population of students.

In the quest to provide approaches to the calculus appropriate for different kinds of student, we should be mindful of the insights of Cauchy. His combination of Euclidean geometry and dynamic movement links naturally to the perceptual idea of dynamic continuity that can now be represented using interactive computer graphics. His ideas of infinitesimal quantities, generated by computable sequences of numbers that become small, uses the same natural mechanism in which such a process is named and then imagined as an arbitrarily small quantity.

Cauchy's ideas of function, continuity, limit and infinitesimal play a significant role in the evolution of modern ideas of the calculus prior to the introduction of the formal limit in mathematical analysis. Using the perceptual idea of local straightness, dynamic variation can be programmed interactively to enable the learner to explore the changing slope of standard functions to predict the formulae for their derivatives in a meaningful human way. Having experienced the derivatives of standard functions perceptually, it becomes more natural to formulate the modern limit concept as a technique to compute the precise symbolic formulae for the derivatives of combinations of such functions.

The ideas of Cauchy are not only appropriate for his own era, they provide a basis for the evolution of all modern approaches to the calculus, be it through the pragmatic use of the Leibniz notation in applications, the use of interactive dynamic computer graphics to visualise ideas and software programs that manipulate symbolism, or the more rigorous developments in standard or non-standard analysis.